\documentclass[english]{smfart}
\usepackage{smfthm}

\usepackage[latin1]{inputenc}
\usepackage[all]{xy}
\usepackage{amssymb} % \varnothing

\usepackage{hyperref}

\newtheorem{lemma}{Lemma}
\newtheorem{theorem}{Theorem}
\newtheorem*{theorem*}{Theorem}
\newtheorem{definition}{Definition}
\newtheorem{proposition}{Proposition}
\newtheorem{corollary}{Corollary}
\newtheorem*{Claim*}{Claim}

\def\C{\mathbb{C}}
\def\D{\mathbb{D}}
\def\H{\mathbb{H}}
\def\N{\mathbb{N}}
\def\R{\mathbb{R}}
\def\S{\mathbb{S}}
\def\T{\mathbb{T}}
\def\Z{\mathbb{Z}}
\def\ds{\displaystyle}
\def\on{\operatorname}
\def\wt{\widetilde}
\def\wh{\widehat}
\def\tend{\longrightarrow}
\def\cal{\mathcal}

\def\Hl{H}

\renewcommand{\epsilon}{\varepsilon}
\renewcommand{\emptyset}{\varnothing}
\renewcommand{\Im}{\on{Im}}

\newcommand{\Lav}{L}
\newcommand{\LavFatou}{L'}
\newcommand{\FillLav}{N}
\newcommand{\FillLavFatou}{N'}
\newcommand{\ins}[1]{\on{int}(#1)}
\newcommand{\setof}[2]{\big\{{#1}\ \big|\ {#2}\big\}}
\newcommand{\hr}{
\bigskip
\begin{center}\rule{3cm}{0.5pt}\end{center}
\bigskip
}
\newcommand{\rk}{\textbf{Remark.}\ }

\begin{document}

\title[The measure of Julia-Lavaurs sets with a virtual Siegel
disk]{The measure of quadratic Julia-Lavaurs sets with a bounded type
virtual Siegel disk is zero}
\author{Arnaud Chéritat}

\maketitle

We shall prove here what was called hypothesis~3 in the author's
thesis \cite{Cthese}, page~159.

\section{Theorem}

\subsection{In short}

Let $p/q$ be an irreducible fraction ($q>0$). Let us define the polynomial
$P(z)= e^{i 2\pi p/q} z + z^2$.
It is well known that the parabolic fixed point $z=0$ is the only non
repelling periodic point of $P$, and moreover that it has precisely
$q$ repelling petals. The action of $P$ on the (repelling) \'Ecalle-Voronin
cylinders is therefore transitive. We will choose one end of these
cylinders and call it $\nu$.
Let $\omega\in\R$ be any bounded type irrational.
The theory of \emph{parabolic enrichment} then associates to
this data a unique Julia-Lavaurs set $\Lav$ with a virtual Siegel disk
with center $\nu$ and virtual multiplier $\exp(i2\pi \omega)$.

\begin{theorem}\label{thm_main}
  \[\on{leb} \Lav =0.\]
\end{theorem}

It shall be noted that we will not make use of the classification of
Fatou components of Julia-Lavaurs sets.

\subsection{More details}

For typesetting system reasons, we will note $\ins{X}$ the interior of
a set $X$.

Let $K$ be the filled-in Julia set, and $J$ the Julia set, of
\[P(z)= e^{i 2\pi p/q} z + z^2.\]
The theory of parabolic enrichment defines Lavaurs maps $g_\sigma :
\ins{K} \to \C$, where $\sigma \in \C$ is a complex parameter.
They commute with $P$. For simplification, the Lavaurs map we will use is
a modification of the usual Lavaurs map, that is described in
\cite{Cgeom}. It is less symmetric (the relation $g_\sigma \circ P = P
\circ g_\sigma$ is slightly modified, but $g_\sigma \circ P^q = P^q
\circ g_\sigma$ still holds; moreover they define the same Julia-Lavaurs
sets). 

The filled-in Julia-Lavaurs set $N$ consists in the complement of preimages of
the basin of infinity by iterates of the Lavaurs maps:
\[\C \setminus \FillLav = \bigcup_{n\in\N} g_\sigma^{-n} (\C\setminus
K)\]
One of many equivalent definitions of the Julia-Lavaurs set is the
boundary (in $\C$) of the filled-in Julia-Lavaurs set:
\[\Lav = \partial \FillLav.\]
This set is compact, completely invariant by $P$ (meaning $P(L)=L$ and
$P^{-1}(L)=L$), and contains $J$. 

These sets depend only on $\sigma \bmod \Z$, as follows for instance
from the property: $g_{\sigma+1} = P^q \circ g_\sigma$ and the fact
$K$ is completely invariant by $P$.

A closely related and more easily understandable function is the
\emph{horn map} $h_\sigma$. Let $T_\sigma$ be the translation by
$\sigma$ in the complex plane.
The Lavaurs and horn maps are related to
each other by the following relations:
\begin{eqnarray*}
 g_\sigma & = &  \psi_+ \circ T_\sigma \circ \phi_\div \\
 h_\sigma & = & T_\sigma \circ \phi_\div \circ \psi_+
\end{eqnarray*}
Where $\phi_\div :\ins{K} \to \C$ and 
$\psi_+ : \C \to \C$ are analytic functions.
This implies for instance that $\psi_+$ is a semi-conjugacy from
$h_\sigma$ to $g_\sigma$: $\psi_+ \circ h_\sigma = g_\sigma \circ \psi_+$.

\hr

These relations are better explained on the following \emph{non
  commuting} diagram:
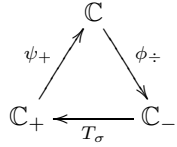
\[
 \xymatrix@=30pt{
  \save[]+<25pt,0pt>*+{\C}="targ" \ar[dr]^{\phi_{\div}} \restore
  & \\
  \C_+ \ar"targ"^{\psi_+} & \ar[l]^{T_\sigma} \C_{-} \\
 }
\]
In this diagram, $\C$ is the complex plane where $P$ lives, $\C_-$ is
the plane of the attracting Fatou coordinates,
$\C_+$ the plane of the repelling Fatou coordinates,
$\phi_\div : \ins{K} \to \C_-$ is the extended attracting Fatou
coordinate, and $\psi_+ : \C_+ \to \C$ is
the extended repelling Fatou coordinate.  See \cite{D}, then
\cite{Cgeom} for their definitions. 
Then $g_\sigma$ and
$h_\sigma$ are defined by following the arrows according to the
illustration below:
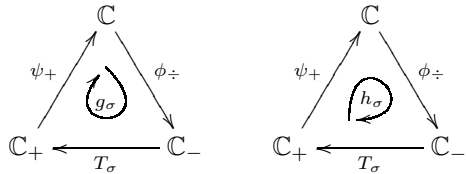
\[
 \qquad
  \xymatrix@=40pt{
  \save[]+<30pt,0pt>*+{\C}="targ" \ar[dr]^{\phi_{\div}} \restore
  \save[]+<30pt,-20pt>*+{}="grat" \ar@(dr,dl)"grat"_{g_\sigma} \restore
  & \\
  \C_+ \ar"targ"^{\psi_+} & \ar[l]^{T_\sigma} \C_{-} \\
 }
 \qquad
 \xymatrix@=40pt{
  \save[]+<30pt,0pt>*+{\C}="targ" \ar[dr]^{\phi_{\div}} \restore
  \save[]+<22pt,-40pt>*+{}="grat" \ar@(u,r)"grat" \restore
  \save[]+<31pt,-32pt>*+{{}_{h_\sigma}}="gart" \restore 
  & \\
  \C_+ \ar"targ"^{\psi_+} & \ar[l]^{T_\sigma} \C_{-} \\
 }
\]

\hr

The preimage of $\Lav$  in repelling Fatou
coordinates is a closed set, $\psi_+^{-1}(\Lav)$, that we will note
$\LavFatou$. The preimage of $\FillLav$ will be noted $\FillLavFatou$.
And we set $J' = \psi_+^{-1}(J)$ and $K' = \psi_+^{-1}(K)$, which are
also closed. Then $L'$ and $N'$ can be characterized in terms of
$h_\sigma$ and $K'$: 

\begin{enumerate}
\item\label{item_e1} $\on{def} h_\sigma=\ins{K'}$,\smallskip 
\item\label{item_e2} $\C \setminus \FillLavFatou = \bigcup_{n\in\N}
h_\sigma^{-1}(\C\setminus K')$,\smallskip
\item\label{item_e3} $\LavFatou = \partial \FillLavFatou$.
\end{enumerate}

\proof
Obviously, $\on{def} h_\sigma = \psi_+^{-1} (\ins{K})$.
Claim (\ref{item_e1}) then follows from $\psi_+$ being a continuous
and open map: $\psi_+^{-1} (\ins{K}) = \ins{\psi_+^{-1}(K)}$.
Claim (\ref{item_e2}) is trivial.
Claim (\ref{item_e3}) also follows from $\psi_+$ being continuous and open.
\qed

\medskip

From the definition of $h_\sigma$, it follows that
\[h_\sigma^{-1}(N') = N' \cap \on{def} h_\sigma.\]
Since $h_\sigma$ is an open and continuous map, this yields
\[h_\sigma^{-1}(L') = L' \cap \on{def} h_\sigma.\]

Let $\pi : \C \to \C/\Z$ be the canonical projection.
The map $h_\sigma$ turns out to commute with $T_1$. It therefore defines
a quotient map from $\pi(\ins{K'})$ to $\C/\Z$.
The cylinder is identified to $\C^*$ via the
exponential map $z \mapsto \exp(i 2\pi z)$. Let $\S^2$ be the Riemann
sphere: the embedding $\C^* \subset \S^2$, enables to consider 
the cylinder $\C/\Z$ completed by its two ends as a Riemann surface
isomorphic to $\S^2$. We will note $+i\infty$ the upper end
(corresponding to $0$ in $\S^2$) and
$-i\infty$ the lower end (corresponding to $\infty$ in $\S^2$).

The set $\ins{K'}$ is a neighborhood of both ends of the cylinder, and
it turns out that $h_\sigma$ has an analytic extension
there, fixing each end, with non zero multiplier (see \cite{DH},
\cite{Lav}, \cite{S}). These two
multipliers are called the \emph{virtual multipliers}. 

Now, remind that the $h_\sigma$ are defined as the following
composition: $T_\sigma \circ \phi_\div \circ \psi_+$. Therefore, given
$t\in\C$, changing $h_\sigma$ into $h_{\sigma + t}$ changes the
virtual multipliers by multiplying them by $\exp(i 2\pi t)$ for the
upper end and $\exp(-i 2\pi t)$ for the lower end. This is why, given
a rotation number $\omega$ and one of the two ends of the cylinder, there
is only one value of $\sigma \bmod 1$ such that $h_\sigma$ has virtual
multiplier $\exp(i 2\pi \omega)$ at this end. From now on, we will fix
$\sigma$ to this value.

\section{Proof}

The question of measure of the Julia set in the case of a quadratic
polynomial with a bounded type Siegel disk has been first settled by
Petersen and Lyubich. They use a quasiconformal model constructed from a
Blaschke fraction. Then McMullen proved this Julia set has Hausdorff
dimension $< 2$. These works have been further generalised.

In our case, points in the Julia set will be separated
into classes, according to their behavior: $P$-type, $J$-type,
$R$-type and $F$-type points. The first class will be handled by adapting
McMullen's techniques to our setting: it turns out to be quite
straightforward. The other classes are easier (in the
author's opinion).

We will make use of the following notations. If $I\subset\R$ is an open
interval, $\C_I$ is defined as $C_I=\C\setminus(\R\setminus I)$. It is
an open subset of $\C$. If $U\subsetneq\C$ is a connected and simply
connected open subset of $\C$, let $d_U$ be the hyperbolic distance on
$U$ and $B_U(z,r)$ the hyperbolic ball of center $z$ and radius $r$.

\subsection{Remark about critical circle maps}

Let us recall the

\begin{definition}
  Here, a \emph{critical circle map} is a map $f : \T \to \T$ that is
  analytic, bijective, orientation preserving, and has one and only
  one critical point, of local degree $3$.
\end{definition}

A reference on this is \cite{dFdM}, that generalise
the work of \cite{McMu}, using \cite{Ya}. 

Here, for simplicity, we chose to adapt \cite{McMu} to prove that the
set of $P$-type points has measure equal to $0$. However,
the techniques of Lyubich explained in \cite{Ya} probably directly
give it. The author ignores if \cite{dFdM} does.

\subsection{Reduction to a model}\label{subsec_redmod}

To fix ideas, we will now assume that the end of the cylinder at which
we put the Siegel disk is the upper end. All the discussion below works
for the other case.

In~\cite{Cgeom}, the author defined a pre-model map $\beta$.
It plays for $h_\sigma$ the role that the Blaschke fraction 
$B=z^2\frac{z-3}{1-3z}$ plays for the quadratic map $e^{i 2\pi \omega}
z + z^2$ in the Douady, Ghys, Herman (and others) surgery: in both
cases one can modify by surgery the pre-model to obtain a model
$\wt{\beta}$ (resp.\ $\wt{B}$). See~\cite{Cgeom} for details. 

We will use the following properties of $\beta$ and $\wt{\beta}$: (the
following list has redundancies)
\begin{enumerate}
\item $\on{def}(\wt{\beta})$ is an open set that contains the closed
  upper half plane $``\Im z \geq 0"$,
\item $\wt{\beta}$ is a continuous open map,
\item $\beta$ is analytic, and real-symmetric (i.e.\ $\on{def}(\beta)$
 is invariant by the reflection $s : z \mapsto \overline{z}$, and
 $\beta$ commutes with $s$),
\item $\R\subset \on{def}(\beta)$,
\item the restrictions of $\wt{\beta}$ and $\beta$ to the
  half-plane $``\Im z \leq 0"$ coincide (i.e.\ they have the same set
  of definition and are equal on this set),
\item in particular, $\wt{\beta}$ is analytic below $\R$,
\item the restriction of $\wt{\beta}$ to $\R$ induces a critical
  circle map, with rotation number $\omega$,
\item $\wt{\beta}$ induces a homeomorphism from $\H$ to $\H$,
\item $\on{def}(\wt{\beta})$ is $T_1$ invariant and $\wt{\beta}$
  commutes with $T_1$, and the same hold for $\beta$,
\item there is a quasiconformal map $S$, commuting with $T_1$, and
  conjugating $\wt{\beta}$ to $h_\sigma$, sending the upper half
  cylinder to the Siegel disk of $h_\sigma$ associated to the upper
  end,
\item $0$ is a critical value of $\wt{\beta}$.
\end{enumerate}

The maps $h_\sigma$ and $\wt{\beta}$ commute with $T_1$ and therefore
induce dynamics on the cylinder. The map $S$ commutes with $T_1$ and
therefore induces a (quasiconformal) conjugacy between them, that fixes
both ends of the cylinder. Let us call $\wh{h}_\sigma$ the
continuation (of the map induced by $h_\sigma$) that fixes both
ends. Thus the map induced by $\wt{\beta}$ has also a continuation
$\wh{\beta}$ that fixes both ends. It is analytic at the lower end.

Since quasiconformal maps preserve sets of zero Lebesgue measure, it
is enough to work on the model $\wt{\beta}$. If we define sets $J''$, $K''$,
$L''$, $N''$ as the preimages of $J'$, $K'$, $L'$, $N'$ by $S$, we have an
analogous characterization in terms of $\wt{\beta}$.

\begin{enumerate}
\item $\on{def} \wt{\beta}=\ins{K''}$,\smallskip
\item $\C \setminus N'' = \bigcup_{n\in\N}
\wt{\beta}^{-1}(\C\setminus K'')$,\smallskip
\item $L'' = \partial N''$.
\end{enumerate}

We will note $b$ be the restriction of $\wt{\beta}$ to the half plane
$H=``\Im z<0"$: it is an analytic map from an open subset of $H$ to $\C$.

\subsection{Covering properties}\label{subsec_covprop}

Let $X$, $Y$ be Riemann surfaces.
For an analytic function $f : X \to Y$,
let $V_a(f)$ be its set of asymptotic 
values, and $V_c(f)$ its set of critical values. We recall that
asymptotic values are points $y \in Y$ for which there exists a
continuous function $\gamma : [0,+\infty[ \to X$ going to
$X$'s infinity (i.e. eventually avoiding all compact) with
$g(\gamma(t)) \underset{n\to +\infty}{\tend} y$.
Note that the set of asymptotic values depends on which set $Y$ is
being considered: replacing $Y$ by a Riemann surface that
contains it happens to increase the set of asymptotic values of $f$.
Taking $Y=\C$, we have
\begin{enumerate}
  \item $V_a(h_\sigma) = \emptyset$,
  \smallskip
  \item $V_c(h_\sigma) = (\sigma + z_0) + \Z$.
\end{enumerate}
and therefore, 
\begin{enumerate}
  \item $V_a(\beta) = \emptyset$,
  \smallskip
  \item $V_c(\beta) = \Z$.
\end{enumerate}

We will use the following well known property:
\begin{proposition}\label{prop_rev}
  Assume $f : X \to Y$ is an analytic map between Riemann surfaces. If
  $V$ is a simply connected connected open subset of $Y$ that contains
  no asymptotic value and no critical value, then for all 
  connected component $U$ of $f^{-1}(V)$, $f : U \to V$ is an
  analytic isomorphism.
\end{proposition}

\subsection{About the cycles of $\wh{\beta}$}

Apart from the upper end of the cylinder, the cycles of $\wh{\beta}$
are all analytic (meaning $\wh{\beta}$ is
analytic in a neighborhood of the cycle) since they are either equal
to the lower end or contained in the half plane $``\Im z <0"$.

\begin{lemma}
  The lower end of the cylinder is repelling for $\wh{\beta}$.
\end{lemma}
\proof[Proof~1]
  It is known that the product of multipliers of both ends of
  $\wh{h}_\sigma$ (which is independent of $\sigma$) has modulus $>1$
  (see for instance~\cite{BE}).
  Therefore, since the upper end is neutral, the lower end is
  repelling. Now, since there is a topological characterization of
  repelling points, this implies $\wh{\beta}$ is also repelling at
  this point.
\qed
\proof[Proof~2]
  We can give a specific proof in our case: Let $U$ be the connected
  component of the preimage of the half plane $H = ``\Im z"<0$ whose
  projection by $\pi : \C \to \C/\Z$ is a neighborhood of the lower
  end (hence $T_1(U)=U$). Then according to
  proposition~\ref{prop_rev}, $\wt{\beta} : U \to L$ is an isomorphism.
  Therefore, $\wh{\beta} : \pi(U)\cup \{-i\infty\} \to \pi(L)\cup
  \{-i\infty\}$ is an analytic isomorphism between two simply
  connected sets, fixing $-i\infty$. The first being a strict subset
  of the second, $-i\infty$ is repelling.
\qed

\medskip
In fact, A.~Epstein adapted Fatou-Shishikura's inequality to a class
of maps including the horn maps $h_{\tau}$, and since they possess only
one critical value, they can have at most one non repelling cycle. 
We will not use this result here, but since in the case of
$\wt{\beta}$ there is a very simple proof, we will mention it here:
\begin{lemma}
  All the cycles of $\wh{\beta}$ that are not the
  upper fixed point are repelling. Therefore the only non repelling
  cycle of $\wh{h}_\sigma$ is the upper end.
\end{lemma}
\begin{proof}
  Let $n$ be the period, and $z\in H$ a point in the cycle.
  As in proof~2, the set $H$ contains no critical value of
  $\wt{\beta}$, and proposition~\ref{prop_rev} implies $\wt{\beta}^n$
  is an isomorphism from a simply connected set $U$ with $z\in
  U\subsetneq H$, to $H$.
\end{proof}

\subsection{Definition of a piece: $P_0$}

Remind that $0$ is a critical value of $\beta$, and therefore not
a critical point (otherwise $\beta$ would have integer rotation number
on $\R$).
Let $m\in\Z$ be such that $\beta(0) \in ]m,m+1[$.
Let $V = \C_{]m,m+1[}$. Let $U$ be the connected component of
$\beta^{-1}(V)$ containing $0$.
By covering properties of $\beta$, the restriction $\beta :
U \to V$ is an analytic isomorphism. Let $\Hl = ``\Im(z)<0"$.
Let $P_0 = \Hl \cap U$. From $0\in U$ we know than
$P_0\not=\emptyset$.
Therefore,
\[\wt{\beta} : P_0 \to \Hl\]
is an isomorphism. 

\begin{lemma}\label{lem_loin}
The set $\wt{\beta}^{-1}(\Hl) \setminus (P_0+\Z)$ is at positive distance to $\R$.
\end{lemma}

\begin{proof}
Recall that the picture is invariant by $T_1$.
Now, $\beta$ being analytic and having local degree $1$ or $3$ at
points in $\R/\Z$, every point in $\R/\Z$ has a neighborhood $W$
on which $W \cap \beta^{-1}(\Hl) \cap \Hl$ is contained in $U \bmod \Z$.
The claim follows by compactness of $\R/\Z$.
\end{proof}

\begin{lemma}\label{lem_hfinite}
  \[\sup_{z\in P_0} |\Im(z)|<+\infty\]
\end{lemma}
\begin{proof}
  (first proof) Since $\wh{\beta}$ is defined at $-i\infty$, it
  implies there is a component $U'$ of $\wt{\beta}^{-1}(H)$ that
  contains a lower half plane $``\Im z<-h_1"$ for some $h_1>0$. It is
  therefore enough to prove that $U' \cap P_0 = \emptyset$. Since
  $P_0$ is also a component of $\wt{\beta}^{-1}(H)$, if it intersected
  $U'$, it would be equal to $U'$. Therefore, $P_0$ would be invariant
  by $T_1$. Studying the situation at the real critical point $z=m+1$
  of $\beta$ shows that there would be points $z_0\in P_0$
  and $z_1\in T_1(P_0)$ such that $\beta(z_0) = \beta(z_1)$, which
  contradicts injectivity of $\beta : P_0 \to H$.\\
  (second proof)
  It is well known that $\pi(\C\setminus K'')$ is an annulus,
  separating the two ends of the cylinder. It does not intersect
  the closure of $\H$. It therefore separates $P_0$ from the lower end
  of the cylinder. The lemma follows.
\end{proof}

\subsection{Classification of points of $L''$}

What can happen to a point $z\in L''$ when it is iterated by $\wt{\beta}$?
By definition, $z$ never gets to $\C \setminus K''$. So either it is
eventually mapped to $\partial K''$, in which case we will say $z$ is
of type $J$ (Julia). Or all its iterates are defined. It also never
gets to $\H$, otherwise it would belong to the interior of $N''$,
and thus not to $L''=\partial N''$.
If it is eventually
mapped to $\R$, then it must remain there, and we will say $z$ is of
type $R$ (real). If it belongs 
infinitely many times to $\Hl \setminus P_0 \bmod \Z$, we will say $z$ is of
type $F$ (far). Otherwise, it eventualy lands \emph{and} remains in
$P_0 \bmod \Z$: it will be called of type $P$.

Therefore, we have decomposed
\[L'' = J_t \sqcup R_t \sqcup F_t \sqcup P_t\] 
($\sqcup$ means a disjoint union), where
$X_t$ means the set of points $z \in L''$ of type $X$.

\medskip

Let us recall that we note
$b$ be the restriction of $\wt{\beta}$ to the half plane
$H=``\Im z<0"$, and that $b$ is analytic on its domain of definition.

\subsection{Type $J$}

The following fact is now well known and is for instance a consequence
of \cite{DU} or \cite{L}
\begin{lemma}
  \[\on{leb} J = 0.\]
\end{lemma}
Since $\psi_+$ is analytic, $J' = \psi_+^{-1}(J)$ has measure
$0$. Since $S$ is quasiconformal, $J''=S^{-1}(J')$ also has measure
$0$. Now, 
\[J_t= \bigcup_{n\in\N} b^{-n} (J'').\]
Since $b$ is analytic,
\begin{proposition} 
  \[\on{leb} J_t = 0.\]
\end{proposition}

\subsection{Type $R$}

Of course, $\R$ has measure $0$. From
\[R_t = \bigcup_{n\in\N} b^{-n} (\R)\]
(with the convention that $b^{0}(\R) = \R$) it follows that
\begin{proposition} 
  \[\on{leb} R_t = 0.\]
\end{proposition}

\subsection{Type $F$}

For rational maps, there is a theorem of Lyubich \cite{L} stating that the set
of points in the Julia set that do not tend to the postcritical set has
Lebesgue measure equal to $0$. Here we will give a simple adaptation
to the specific case we are studying.

To see the relation with $F_t$, note that the postcritical set of
$\wt{\beta} \bmod \Z$ is $\R/\Z$, and that there exists $\epsilon>0$
such that points of type $F$ must 
infinitely often visit $\Hl_\epsilon = ``\Im z \leq -\epsilon"$ (as
follows from lemma~\ref{lem_loin}).

Let $\Hl=``\Im z<0"$ and let us call $U$ the connected component of
$\wt{\beta}^{-1}(\Hl)$ that is a neighborhood of $-i\infty$ (when
projected to $\C/\Z$). Let $g$ be
the restriction of $\wt{\beta}$ to $U$.

Let $z\in F_t$ and let us prove that $z$ is not a Lebesgue density
point of $F_t$, by providing arbitrarily small balls 
that do not intersect $F_t$, and whose radii are
commensurable to the distance of their centers to $z$. They are
provided by the classical method: taking univalent branches and applying
Koebe's distortion theorem.

Let us introduce the \emph{nest} of $z$: this is the sequence of
\emph{pieces} $P_n(z)$ where $P_n(z)$ is the connected component
containing $z$ of $\wt{\beta}^{-n}(\Hl) = b^{-n}(\Hl)$. This sequence
is decreasing for inclusion. According to proposition~\ref{prop_rev},
the sets $P_n(z)$ are open and simply connected and $b^n : P_n(z) \to H
$ is an analytic isomorphism. 

\begin{lemma}\label{lem_stz}
  For all $z\in L''$ not of type $J$, the distance from $z$ to the
  boundary of $P_n(z)$ tends to $0$ when $n\tend+\infty$.
\end{lemma}
\proof Otherwise, the intersection $A$ of the sets $P_n(z)$ would be a
neighborhood of $z$, included in $N''$ (since all points $z'\in A$
has infinite orbit that never goes to the upper half plane).
Therefore $z$ would belong to $\ins{N''}$ and thus not to $L'' =
\partial N''$, which leads to a contradiction.\qed

\medskip

\begin{Claim*}
  Every point of $U$ is eventually mapped out of $U$ under
  iteration of $b$.
\end{Claim*}
\proof We defined $g$ as the restriction of $\wt{\beta}$ to $U$.
The univalent map $g^{-1} : \Hl \to U$ commutes with $T_1$, and
its quotient, completed by the lower end of the cylinder, is
conjugated by $z \mapsto \exp(- 2i\pi z)$ to an analytic isomorphism
$f$ from the unit disk $\D$ to an open subset $W$ with compact closure
in $\D$, with $f(0)=0$. In that case, it is well known that the
intersection of the images $f^n(\D)$ is reduced to $\{0\}$. Therefore,
$\ds \bigcap_{n\in\N} g^{-n} (\H) = \emptyset$.
\qed

\medskip

Let now $z\in F_t$.
The set $\C\setminus K''$ is open an contained in $\Hl$. Its preimages
by $b^n$ are of course disjoint from $F_t$. Therefore, by
Koebe's distortion theorem and lemma~\ref{lem_stz}, it is enough to
find some $M$ such that, infinitely many times, the orbit of $z$
passes at hyperbolic distance to $\C\setminus K''$ in $\Hl$ less than
$M$:
\[(\forall z\in F_t)\ (\exists M>0)\ (\forall N\in\N)\ (\exists n \geq
N)\ \on{dist}_\Hl(b^n(z), \C\setminus K'')< M.\]

The constant $M$ will here be independent of $z$ and be equal to
the supremum of the hyperbolic distance in $\Hl$ to $\C\setminus K''$,
of points in the set
\[B = \Hl_\epsilon \setminus g^{-1}(\Hl_\epsilon),\]
where $\Hl_\epsilon = \setof{z \in \C}{\Im z \leq -\epsilon}$.
This supremum is finite, because $\C\setminus K''$ is
$T_1$-invariant and $B$ has imaginary part bounded away from $0$ and
$-\infty$.

Now, let $m\geq N$ such that $b^m(z) \in \Hl_\epsilon$. By the
claim\footnote{The claim reads $\bigcap_{n\in\N} g^{-n}(\Hl) =
\emptyset$, which implies $\bigcap_{n\in\N} g^{-n}(\Hl_\epsilon) =
\emptyset$. Therefore $\Hl_\epsilon = \bigcup_{n\in\N} g^{-n}(B)$.
}, there is some $n \geq m$ such that
$b^n(z) \in B$. This yields:
\begin{proposition} 
  \[\on{leb} F_t = 0.\]
\end{proposition}

\medskip
\noindent\rk Since this works for all $\epsilon>0$, the measure of the set
of points which do not tend to $\R$ is equal to $0$: this is
Lyubich's theorem in our particular case.

\subsection{Type $P$}

Let
\[C = \setof{z\in P_0 \bmod \Z}{\text{all iterates of } z \text{ belong to }
P_0 \bmod \Z}\]
then, by definition
\[P_t = L'' \cap \bigcup_{n\in\N} b^{-n} (C).\]
From $\wt{\beta}^{-1}(L'') = L'' \cap \on{def}{\wt{\beta}}$, it follows that
\[ b^{-1}(L'') = L'' \cap \on{def} b.\]
Therefore\footnote{In fact, more is true (even if we will not use it):
Lavaurs' versions of Sullivan's non-wandering theorem and of the
classification of components, applied
to $h_\sigma$, imply $C \subset L''$ (an therefore $\ins{C} = \emptyset$).
Indeed, by the absence of non-repelling cycles of $\wh{h}$
except from the upper end, every component of $\ins{N'}$ is
eventually mapped to the Siegel disk at the upper end. By the
quasiconformal transformation $S$, this turns into an analogous
statement for $\wh{\beta}$.}, 
\[P_t = \bigcup_{n\in\N}L'' \cap  b^{-n} (C) = \bigcup_{n\in\N}b^{-n} (L'' 
\cap C).\]
Let us postpone to the next section the proof of
\begin{lemma}\label{lem_CparC}
  \[\on{leb} L'' \cap C = 0,\]
\end{lemma}
\noindent from which follows that
\begin{proposition} 
  \[\on{leb} P_t = 0,\]
\end{proposition}
\noindent which will end the proof of theorem~\ref{thm_main}.

\subsection{Proof of lemma~\ref{lem_CparC}}
We will endow the lower half-plane $\Hl$ with the hyperbolic metrics
given by
\[\frac{|dz|}{|\Im z|}.\]

Let us state the main lemma
\begin{lemma}\label{lem_main}
There exists $M>1$ such that $\forall x\in\R$, $\forall s \in
]0,1[$, there is a euclidean ball $B=B(z,r)$ with
\begin{enumerate}
 \item $B \subset H$,
 \item $B$ eventually falls in $\H$ under iteration of $b$, and thus
 $B \cap C = \emptyset$,
 \item $r> s/M$,
 \item $|z-x| < s M$.
\end{enumerate}
\end{lemma}

%\rk One usually and informally says that for all $x\in\R$,
%there is a Euclidean ball at every scale $s<1$.

Asumme we have fixed some height $h_0>0$.
In terms of hyperbolic metrics in $H$, lemma~\ref{lem_main} implies
that for every point $z'\in \Hl$ with $\Im z'> - h_0$,
there is a hyperbolic ball $B'$ which is eventually mapped to $\H$
under iteration of $b$, with hyperbolic diameter $C_2$ and
contained in the hyperbolic ball of center $z'$ and radius $C_3$. The
constant $C_2$ depends on $M$. The constant $C_3$ depends on $M$ and $h_0$.

Now let $h_0 = \sup_{z\in P_0} |\Im(z)|$. By lemma~\ref{lem_hfinite},
$h_0<+\infty$. 
For $z\in C$, let $r_n$ be
the distance from $z$ to the boundary of the piece $P_n(z)$. 
For a given $n$, set $z'=b^n(z)$ and let $B'$ be as above.
If we pull back $B'$ by the branch
of $b^{-n}$ mapping $\Hl$ to $P_n(z)$, then according to Koebe's
distortion theorem, $B'$ is mapped to a set which contains a euclidean
ball $B''$ with center within distance $< M' r_n$ to $z$ and diameter
$> r_n/M'$. The constant $M'$ depends on $C_2$, $C_3$, but not on $n$, $z$.
Note that these balls $B''$ are eventually mapped to $\H$ under
iteration of $b$, and thus do not intersect $C$. Now if $z \in C \cap
L''$, then $r_n \tend 0$ (by lemma~\ref{lem_stz}).
The balls $B''$ then prevent $z$ from being a density point of $C$.
This proves lemma~\ref{lem_CparC}.

\subsection{Proof of lemma \ref{lem_main}}

For an interval $I$, $|I|$ will denote its length. Intevals $I$ and
$J$ are called $K$-commensurable if and only if the quotient of their
lengths lies between $1/K$ and $K$.

Let us introduce the \emph{dynamical partition}: let $f$ be a critical
circle map with irrationnal rotation number $\theta$ and let $c \in
\R/\Z$ be its critical point. Let $p_n/q_n$ be the convergents of $\theta$.
The partition of $\R/\Z$ into intervals separated by the points
$f^{-k}(c)$, $0\leq k < q_n+q_{n+1}$ is called the \emph{dynamical
partition at level $n$}. Let us fix by convention that the intervals
are of the form $[a,b[$. For $x\in\R/\Z$, we will note $I_n(x)$ the
unique element of the dynamical partition at level $n$ that contains $x$.

Le us recall the following (see~\cite{Yo})
\begin{theorem*}[Yoccoz]
  The map $f : \R/Z \to \R/\Z$ is conjugated to the rotation
  $T_\theta : \R/Z \to \R/\Z$ (by an orientation preserving homeomorphism of $\R/\Z$).
\end{theorem*}
\noindent In particular, the backward orbit $f^{-k}(c)$ is dense in $\R/\Z$.

We will make use of the following theorem of Herman and Swiatek:
\begin{theorem*}[Herman, Swiatek, real bounds]
For all critical circle map $f$ with irrationnal rotation number,
there exists $K>1$ such that for all $n\in\N$, 
any two consecutive intervals of the dynamical partition at level
$n$ are $K$-commensurable.
\end{theorem*}

A well known corollary is the following:
\begin{corollary}\label{cor_part}
If, furthermore, the rotation number has bounded type, then there
exists $K'>1$ such that $\forall x\in\R$, $\forall \ell \in ]0,1[$, $\exists
n\in\N$ such that $|I_n(x)|$ is $K'$-commensurable to $\ell$.
\end{corollary}

The following lemma is purely geometric
\begin{lemma}\label{lem_cone}
  For all $K>1$ there exists $M_0>0$ and $r_0>0$ such that the
  following holds. Assume we are given three intervals $[a,b[$,
  $[b,c[$ and $[c,d[$ in $\R$, and that touching intervals are
  $K$-commensurable. Let $U = \C_{]a,d[}$. Then, \emph{for the
  hyperbolic metrics on $U$}, every cone based on a point
  $x\in[a,d]$, with central direction $-90°$ and opening
  $30°$, contains a hyperbolic ball $B=B_U(z,r)$,
  with radius $r=r_0$, and such that the hyperbolic diameter
  $\on{diam}_U([b,c]\cup B)<M_0$.
\end{lemma}
The proof is elementary. This would hold if the opening angle is
replaced by any number between $0°$ and $180°$.

\begin{lemma}\label{lem_pascone}
  There is some $h_1>0$ such that the set $P_0+\Z$ avoids
  the interior of the triangle of vertices $c$, $c-a-ih_1$, and
  $c+a-ih_1$, where $c\in\R$ is a real critical point of $\beta$, and
  $a$ is chosen so that the angle at vertex $c$ is equal to $30°$.
\end{lemma}
\begin{proof}
  This immediately follows from the critical point of $\beta$ having
  local degree $3$ and $\beta$ sending $\R$ to $\R$ increasingly.
\end{proof}

\begin{corollary}\label{cor_comm}
  For the map $\beta$, there exists $M_2>0$ and $r_2>0$ such that $\forall
  n \geq 2$, for all intervals $I$ of the dynamical partition at
  level $n$, there is a euclidean ball $B=B(z,r)$ eventually mapped to
  $\H$ under iteration of $b$, with
  $r > r_2 |I|$ and $d(z,I)<M_2|I|$.
\end{corollary}
\begin{proof}
  Let us note $I=[b,c[$ and let $[a,b[$ be the previous interval and
  $[c,d[$ the next one, in the dynamical partition at level $n$.  Let
  us note $a=\beta^{-m_1}(u)$, $b=\beta^{-m_2}(u)$,
  $c=\beta^{-m_3}(u)$ and $d=\beta^{-m_4}(u)$, where $u$ is the
  critical point of $\beta$ on $\R/\Z$.  We took $n\geq 2$, therefore
  $q_n+q_{n+1} \geq 5$, so there are at least $4$ distinct points
  defining the dynamical partition at level $n$,
  so the points $a$, $b$, $c$, $d$ are distinct. Thus the $m_j$
  are distinct. Let $m\in\N$ be the least\footnote{in fact, we could
  also take $m = \min(m_2,m_3)$} of $m_1$, $m_2$, $m_3$ and $m_4$. Let
  $[a',b'[$, $[b',c'[$ and $[c',d'[$ be the image of the three
  intervals by $\beta^m$. These are still three consecutive intervals
  of the dynamical partition at level $n$. (Indeed, the image by
  $\beta$ of an interval $I$ in the partition ${\cal P}$ created by
  the preimages of the critcal point, from order $0$ up to $k-1$, is
  still in ${\cal P}$, unless the $k$-th preimage of the critical
  point belongs to $\ins{I}$, which for $k=q_n+q_{n+1}$ implies that
  the critical point bounds $I$.)
%%   NOTE QUE TOUT CE QU'ON A BESOIN C'EST QUE CES 4
%%   NOUVEAUX SONT $K'$-commensurables POUR UN AUTRE $K'$, et ça c'est
%%   automatique en type borné, car on acquiert au plus un mauvais point.
  Therefore, $[a',b'[$ and $[c',d'[$ are $K$-commensurable to
  $[b',c'[$, where $K$ is given by the real bounds.
  Let us consider, for the hyperbolic metrics on $U = \C_{]a',d'[}$,
  the hyperbolic ball $B$ provided by 
  lemma~\ref{lem_cone} for the cone based on the critical point. Now,
  either $B$ is contained in the triangle provided by
  lemma~\ref{lem_pascone}, and we set $B'=B$. If it is not, then $n$
  must be less than some $n_0\in\N$ (indeed, the size of the three
  intervals is bounded from above by a sequence that tends to $0$ as
  $n$ tends to infinity). We then take $B'=$ the biggest hyperbolic disk
  contained in the triangle.
  Then, there exist constants $M_1$ and $r_1$ independant of $n$ such
  that, in both cases, $B'$ has hyperbolic radius (in $U$) $\geq r_1$ and
  $\on{diam}_U(B'\cup[b',c']) < M_1$. Let us now apply the
  branch $g$ of $\beta^{-m}$ defined on $U$ and sending $]a',d'[$ to
  $]a,d[$. This branch exists because $U$ does not contain any
  singular value of $\beta^m$ (otherwise, $]a,d[$ would contain a
  point of the form $f^{-k}(c)$ for some $0\leq k <m$, which is not
  the case).
  The lemma follows using bounded distortion of $g$ in $B_U(b',M_1)$,
  which contains both $B'$ and $[b',c']$.
\end{proof}

Corollaries~\ref{cor_part} and~\ref{cor_comm} give lemma~\ref{lem_main}.

\section*{Acknowledgements}

I would like to thank M.\ Yampolsky for useful discussions.

\end{document}